%Paper: alg-geom/9208002
%From: ribet@math.berkeley.edu (Kenneth A. Ribet)
%Date: Thu, 6 Aug 92 15:54:15 PDT

% This document was prepared for AMS-TeX 2.1 and uses the 2.1 version
% of amsppt.sty.  In fact, I seem to be using version 2.1a of amsppt.sty,
% which is dated 22 Jan 1992.

 \input amstex %in case we're using plain TeX
 \documentstyle{amsppt}

 \def\startproof{\noindent{\sl Proof.}\enspace}
 \def\endproof{{\unskip\nobreak\hfil\penalty50\hskip1em
   \hbox{}\nobreak\hfil$\blacksquare$%
     \parfillskip=0pt \finalhyphendemerits=0 \par}}

 \def\leqno(#1){\eqno{\hbox to .01pt{}\rlap{\rm \hskip -\displaywidth(#1)}}}

 %  This makes the paper look better in my humble opinion
           \parskip\smallskipamount

 % Numbering references

 \newcount\refCount
 \def\newref#1 {\advance\refCount by 1
 \expandafter\edef\csname#1\endcsname{\the\refCount}}

   \newref CHI %
   \newref DELCORVALIS
   \newref FALT
   \newref FW
   \newref BHG
   \newref GRBU
   \newref SGA
   \newref HIDA
   \newref GAD
   \newref LORENZ
   \newref MILNE
   \newref MUMF
   \newref OGG
   \newref ANNALS
   \newref MYTHESIS
   \newref BONN
   \newref MATHANN
   \newref DPP
   \newref MG
   \newref SERREINV
   \newref JPSDURHAM
   \newref DUKE
   \newref DISCONT
   \newref PINK
   \newref SHIM
   \newref CFHO
   \newref FACTOR
   \newref SHTAN
   \newref WEILFD
   \newref WEIL

 % math macros

  \let\tempcirc=\circ
  \def\circ{\mathord{\raise0.25ex\hbox{$\scriptscriptstyle\tempcirc$}}}

  \def\seteq{\mathrel{:=}}
  \def\GalQ{\Gal(\Qbar/\Q)}
  
  \def\Aut{\mathop{\roman{Aut}}\nolimits}
  \def\deg{\mathop{\roman{deg}}\nolimits}
  \def\Br{\mathop{\roman{Br}}\nolimits}

  \def\Gal{\mathop{\roman{Gal}}\nolimits}
  \def\Res{\mathop{\roman{Res}}\nolimits}
  \def\End{\mathop{\roman{End}}\nolimits}
  \def\Lie{\mathop{\roman{Lie}}\nolimits}
  \def\Hom{\mathop{\roman{Hom}}\nolimits}

  \def\GL{\mathop{\bold{GL}}\nolimits}
  
  \def\Norm{\mathop{\bold N}\nolimits}
  \def\GLTWO{$\GL_2$}

  \def\trace{\mathop{\roman{tr}}\nolimits}

  \def\Q{{\bold Q}}

  \def\Qell{\Q_\ell}
  \def\Qellbar{\overline{\Q}_\ell}
  \def\Zell{\Z_\ell}
  \def\C{{\bold C}}
  \def\Qbar{\overline\Q}

  \def\F{{\bold F}}

  \def\Fell{{\F_\ell}}
  \def\Zell{\Z_\ell}

  \def\rholambdabar{\overline{\rho}_\lambda}
  \def\ssigma{{}^\sigma\mkern-\thinmuskip}
  \def\ttau{{}^\tau\mkern-\thinmuskip}
  \def\ggg{{}^g\mkern-\thinmuskip}
  \def\Flambdabar{\overline{\F}_\lambda}
  \def\chiellbar{\overline{\chi}_\ell}
  \def\sigmabar{{\bar{\sigma}}}
  \def\Z{{\bold Z}}
  
  \def\Frob{{\roman{Frob}}}
  \def\Xo#1{X_o(#1)}  % ditto
    % ditto
  \let\phi\varphi

 \topmatter

  \title
 Abelian varieties over $\Q$ and modular forms
  \endtitle

  \author
 Kenneth A. Ribet
  \endauthor

  \affil
 University of California, Berkeley
  \endaffil
 \rightheadtext{Abelian varieties over $\Q$ and modular forms}
 \address{UC Mathematics Department, Berkeley, CA  94720 USA}
 \endaddress

 \email ribet\@math.berkeley.edu \endemail

   \thanks
 This manuscript was written in conjunction with the author's
 lectures at the seventh KAIST mathematics workshop in
 Taejon, Korea (August 11--14, 1992).  It exposes work which
 was done during the author's visits to the Hebrew University
 of Jerusalem and the Universit\'e de Paris~XI (Orsay).  It
 is a pleasure to thank these institutions for their
 hospitality.
 The author wishes to thank J-P.~Serre for suggesting the
 problem and more pointing out the relevance of Tate's
 theorem (Theorem~6.3 below.)
 The author was supported in part by NSF Grant \#DMS 88-06815.
   \endthanks
 \endtopmatter

 \document

 \head 1. Introduction.\endhead \noindent
 Let $C$ be an elliptic curve over~$\Q$.
 Let $N$ be the conductor of~$C$.
 The Taniyama conjecture
 asserts that there is a non-constant map of algebraic curves
 $\Xo N \to C$ which is defined over~$\Q$.
 Here, $\Xo N$
 is the standard modular curve associated with the problem
 of classifying elliptic curves $E$ together
 with cyclic subgroups of~$E$ having
 order~$N$.

 The Taniyama  conjecture may be reformulated in various
 ways.  For example, let $X_1(M)$ be the modular curve
 associated with the problem of classifying elliptic curves
 $E$ along with a point of order~$M$ on~$E$.
 One knows that
 if there is a non-constant map
 $X_1(M) \to C$ over~$\Q$ for some~$M\ge1$, then there is a
 non-constant map $\Xo N\to C$.  (For a result in this
 direction, see \cite\LORENZ.)

 In a recent article \cite\GAD, B.~Mazur introduced another
 type of reformulation.
 Mazur says that $C$ possesses a ``hyperbolic uniformization
 of arithmetic type" if there is a non-constant map
 over the complex field $\pi\:X_1(M)_{\C}\to C_{\C}$ for some
 $M\ge1$.  The map~$\pi$
 is not required to be defined over~$\Q$.
 In \cite{\GAD, Appendix},
 Mazur proves: Suppose that
 there is a non-constant map $\pi\:X_1(M)_{\C}\to C_{\C}$.
 Then there is a non-constant map $\pi'\:X_1(M') \to C$
 over~$\Q$, where $M'$ is a suitable positive integer
 (which may be different from~$M$).
 As mentioned above,
 it follows easily from
 the existence of~$\pi'$ that $C$ satisfies Taniyama's conjecture.
 We thus arrive at the following conclusion:

 \proclaim{(1.1) Theorem {\rm [Mazur]}}
 Let $C$ be an elliptic curve over~$\Q$.  Then $C$
 satisfies Taniyama's conjecture if and only if
 $C$ admits
 a hyperbolic uniformization of arithmetic type over~$\C$.
 \endproclaim

 In connection with \cite\GAD, Serre asked for a conjectural
 characterization of elliptic curves over $\C$ which
 admit a hyperbolic uniformization of arithmetic type.  All
 such curves are defined over~$\Qbar$, the algebraic closure
 of~$\Q$ in~$\C$.  Thus one wishes to characterize
 those elliptic curves over~$\Qbar$ which are quotients of some
 modular curve $X_1(M)$.
 Equivalently, one wishes to study the set of elliptic
 curves
 which are quotients of the Jacobian $J_1(N)_{\Qbar}$
 of~$X_1(N)_{\Qbar}$, for some $N\ge1$.
 It is well known that all complex multiplication elliptic
 curves have this property; this follows from \cite{\PINK,
 Th.~1}.
 Hence we are principally interested in elliptic curves
 $C/\Qbar$ without complex multiplication.

 In this article, we introduce the concept of an
 abelian variety over~$\Q$ which is of ``\GLTWO-type."
 Roughly speaking, this is an abelian variety over~$\Q$
 whose algebra of $\Q$-endomorphisms is a
 number field of degree equal to the dimension of the
 abelian variety.  It is easy to see that $J_1(N)$ decomposes
 up to isogeny over~$\Q$ as a product of such abelian
 varieties.  Hence any $C/\Qbar$ which is a quotient of some
 $J_1(N)$ is a quotient (over~$\Qbar$)
 of an abelian variety of \GLTWO-type over~$\Q$.

 We prove in this
 article the following facts (the third is a
 simple corollary of the first two):
 \item{1.} Let $C$ be an elliptic curve over~$\Qbar$
 without complex multiplication.  Then
 $C$ is a quotient of an abelian variety of~\GLTWO-type over
 $\Q$ if and only if $C$ is isogenous to each of its Galois
 conjugates $\ssigma C$ with $\sigma\in\GalQ$.
 (If $C$ has this latter property, we say that $C$
 is a ``$\Q$-curve."  The terminology is borrowed from
 Gross \cite\BHG.)
 \item{2.} Assume Serre's
 conjecture \cite{\DUKE, (3.2.4$_?$)}
 on representations of~$\GalQ$.
 Then every abelian variety of \GLTWO-type over~$\Q$ is a
 quotient of~$J_1(N)$ for some $N\ge1$.
 \item{3.} Let $C/\Qbar$ be an elliptic curve without complex
 multiplication.  If $C$ is a quotient of $J_1(N)_{\Qbar}$ for
 some $N$, then $C$ is a $\Q$-curve.
 Conversely, suppose that $C$ is a $\Q$-curve.  Then
 if the conjecture \cite{\DUKE, (3.2.4$_?$)} is correct, $C$
 is a quotient of~$J_1(N)_{\Qbar}$ for some~$N$.

 \noindent
 In summary, we arrive at a conjectural characterization of
 elliptic curves over~$\Qbar$ which are quotients of
 modular curves $X_1(N)$: they are exactly the
 $\Q$-curves in the sense indicated above.
 This characterization was predicted by~Serre.
 \head 2.
 Abelian varieties over $\Q$ of \GLTWO-type.\endhead \noindent
 We will be concerned with abelian varieties over~$\Q$ which
 admit actions of number fields that are ``as large as possible."
 We first quantify this concept.

 Suppose that $A$ is an abelian variety over~$\Q$ and that $E$
 is a number field acting on $A$ up to isogeny
 over~$\Q$:
 $$ E \hookrightarrow \Q \otimes_\Z \End_\Q(A).$$
 By functoriality, $E$ acts on the space of tangent vectors
 $\Lie(A/\Q)$, which is a $\Q$-vector
 space of dimension $\dim A$.
 (For an account of the Lie algebra attached to a group
 scheme over a field, see for example \cite{\MUMF, \S11}.)
 The dimension of this vector space
 is therefore a multiple of~$[E\:\Q]$, so that we have
 $[E\:\Q] \mid \dim A$.  In particular, $[E\:\Q] \le \dim A$.

 This observation motivates the study of abelian varieties $A/\Q$
 whose endomorphism algebras contain number fields of maximal
 dimension $\dim A$.  If $E$ is a number field of degree~$\dim A$
 which is contained in~$\Q\otimes\End_\Q(A)$,
 then the Tate modules $V_\ell(A)$ associated
 with~$A$ are free of rank~two over $E\otimes_\Q\Qell$.  Accordingly,
 the action of~$\GalQ$ on~$V_\ell(A)$ defines a representation
 with values in~$\GL(2,E\otimes\Qell)$.  We
 say that $A$ is of ``\GLTWO-type."

 Suppose that $B$ is of \GLTWO-type, and that
 $F \subseteq \Q\otimes\End_\Q(B)$ is a number field of
 degree $\dim B$.  Let $E$ be a number field containing~$F$,
 and let $n=[E\:F]$.  After choosing a basis for $E$ over~$F$,
 we find an embedding $E \subseteq M(n,F)$ of~$E$ into the
 ring of $n$ by~$n$ matrices over~$F$.  Since $M(n,F)$ acts
 naturally up to isogeny
 on $A\seteq B\times\cdots\times B$ ($n$ factors), we obtain an
 embedding $E\subseteq \Q\otimes\End_\Q(A)$.  Hence $A$ is again
 of~\GLTWO-type. (We can summarize the situation by writing
 $A=E\otimes_F B$.)

 We say that an abelian variety $A/\Q$ of \GLTWO-type is
 {\it primitive\/} if it is not isogenous over~$\Q$ to
 an abelian variety
 obtained by this matrix construction, relative to an extension
 $E/F$ of degree~$n>1$ (cf.~\cite{\SHTAN, \S8.2}).

 \proclaim{(2.1) Theorem}
 Let $A$ be an abelian variety of
 \GLTWO-type over~$\Q$.
 Then the following conditions are equivalent: (i) $A$
 is primitive; (ii) $A/\Q$ is simple; (iii) the endomorphism
 algebra of~$A/\Q$ is
 a number field whose degree
 coincides with the dimension of~$A$.
 \endproclaim
 \startproof
 Let $E$ be a number field of degree~$\dim A$ which is
 contained in the $\Q$-algebra
 $\Cal X\seteq\Q\otimes_\Z \End_\Q(A)$.
 Let $D$ be the commutant of~$E$ in~$\Cal X$.
 We claim that $D$ is a division algebra (cf.~\cite{\ANNALS,
 Th.~2.3}).

      To prove the claim,
 we must show that each non-zero
 $\Q$-endomorphism of~$A$ which
 commutes with~$E$ is an isogeny.
 Let $\lambda$ be such an endomorphism and
 let $B$ be the image
 of~$\lambda$. Then $B$ is a non-zero abelian subvariety of~$A$.
 The field $E$ operates on $B$ (up to isogeny), and thereby
 acts by functoriality on the $\Q$-vector space $\Lie(B/\Q)$.
 The dimension of $\Lie(B/\Q)$ is accordingly
 a multiple of~$[E\:\Q]$; on the other hand, it
 coincides with~$\dim B$.  Hence $B=A$, so that $\lambda$ is
 an isogeny.

 The Lie algebra $\Lie(A/\Q)$ may now be viewed as
 a $D$-vector space. Because of this, the dimension of~$\Lie(A/\Q)$
 is a multiple of
 $\dim_\Q(D)$.
 In other words, we have $\dim(D) \mid \dim(E)$.
 This gives the equality $D=E$; i.e., it shows that $E$ is its own
 commutant in~$\Cal X$.

 In particular, the center $F$ of~$\Cal X$
 is a subfield
 of~$E$.  We have then $\Cal X \approx M(n,Q)$,
 where $Q$ is a division algebra with center~$F$.  If $Q$ has dimension
 $t^2$ over~$F$, then $nt = [E\:F]$.  This follows from the fact
 that $E$ is a maximal commutative semisimple subalgebra
 of~$\Cal X$,
 which in turn follows from
 the statement that $E$ is its own commutant in~$\Cal X$.

 The structure of~$\Cal X$
 shows that~$A$ is
 isogenous (over~$\Q$) to a product of $n$ copies of an abelian
 variety~$B$ whose algebra of $\Q$-endomorphisms is
 isomorphic to~$Q$.
 The same Lie algebra argument we have already used shows
 that $\dim_\Q(Q) \mid \dim(B)$, so that $n\cdot\dim_\Q(Q) \mid \dim(A)$.
 This gives the divisibility $n t^2 [F\:\Q] \mid [E\:\Q]$, which
 implies
 $nt^2 \mid nt$.  One deduces that $t=1$, i.e., that
 $Q=F$, and obtains the equality $n=[E\:F]$.  Hence the dimension of~$B$
 coincides with the degree of $F$ over~$\Q$.  Also, we have
 $\Q\otimes\End_\Q(A) \,\approx\, M(n,F)$.

 In other words, $A$ is obtained from $B$ and~$F$ by the construction
 we outlined above.  The equivalence of the three statements in
 the theorem is now clear: each assertion is equivalent to the
 equality $n=1$.
 \endproof

 \head 3. $\ell$-adic representations attached to primitive
 abelian varieties over~$\Q$
 of \GLTWO-type. \endhead\noindent
 In what follows, we study primitive
 abelian varieties of \GLTWO-type over~$\Q$.
 Since we never encounter such abelian varieties which
 are not primitive in the above sense, we will often
 drop the word ``primitive" and refer simply to
 abelian varieties over~$\Q$ of \GLTWO-type.

    To motivate the study of these varieties, we first
 allude to the existence of a large
 class of examples of (primitive) abelian varieties of
 \GLTWO-type over~$\Q$.
 Suppose that $f=\sum a_n q^n$ is
 a normalized cuspidal
 eigenform of weight two on a subgroup of
 $\bold{SL}(2,\Z)$ of the form $\Gamma_1(N)$.
 Then Shimura~\cite{\SHIM, Th.~7.14}
 associates to~$f$ an abelian variety $A=A_f$
 over~$\Q$ together with an action on~$A$ of the field
 $E=\Q(\ldots,a_n,\ldots)$.  The variety $A_f$ may be
 constructed as a quotient
 of $J_1(N)$, the Jacobian of the standard
 modular curve $X_1(N)$ \cite{\FACTOR}.
 The dimension of~$A$ and the degree of~$E$ are equal.
 It is well known (and easy to show)
 that $E$ is the full algebra of endomorphisms
 of~$A$ which are defined over~$\Q$ \cite{\MATHANN, Cor.~4.2}.
 Thus, $A$ is of \GLTWO-type over~$\Q$.
 For each $N\ge1$, the Jacobian $J_1(N)$
 is isogenous to a product of abelian varieties of the
 form~$A_f$ \cite{\MATHANN, Prop.~2.3}.

 In \S4, we
 study the converse problem: Suppose that
 $A/\Q$ is of \GLTWO-type.
 Is $A$ isogenous to
 a quotient of~$J_1(N)$, for some $N\ge1$?
 We show that an affirmative answer is implied by
 conjecture~(3.2.4$_?$) of Serre's article
 \cite\DUKE\ on modular representations of~$\GalQ$.

 We now begin the study of Galois representations
 attached to
 \GLTWO-type abelian varieties over~$\Q$.
 Suppose that $A$ is such an abelian variety.
 Let $E$ be the endomorphism algebra of~$A/\Q$.
 Then $E$ is a number field which is
 either a totally real number field
 or a ``CM field," since each $\Q$-polarization of~$A$
 defines a positive involution on~$E$.

 Recall that
 for each prime number $\ell$, the Tate module $V_\ell=V_\ell(A)$
 is free of rank two over $E\otimes_\Q \Qell$.
 For each prime $\lambda \mid\ell$ of~$E$, let $E_\lambda$
 be the completion of~$E$ at~$\lambda$, and set
 $V_\lambda\seteq V_\ell\otimes_{E\otimes\Qell} E_\lambda$.
 Thus $V_\lambda$ is a two-dimensional vector space over~$E_\lambda$,
 and $V_\ell$ is the direct sum of the~$V_\lambda$ with $\lambda\mid\ell$.
 The action of~$\GalQ$ on~$V_\lambda$ defines a ``$\lambda$-adic
 representation"~$\rho_\lambda$.  One knows that the collection
 $(\rho_\lambda)$ (as $\lambda$ ranges over the set of finite
 primes of~$E$) forms a strictly compatible system of $E$-rational
 representations whose exceptional set is the set of prime
 numbers at which $A$ has bad reduction.
 (For background on this material,
 see \cite{\DISCONT, \S11.10} and perhaps
 \cite{\MYTHESIS, Ch.~II}.)
 We will prove some facts about the $\lambda$-adic representations
 $\rho_\lambda$ and their reductions mod~$\lambda$.

 For each $\lambda$, let $\delta_\lambda\:\GalQ\to E^*_\lambda$
 be the determinant of~$\rho_\lambda$.
 The $\delta_\lambda$ form a compatible
 system of $E$-rational one-dimensional representations
 of~$\GalQ$.  In the case where $E$ is totally real, we have
 $\det\rho_\lambda=\chi_\ell$, where $$\chi_\ell\:\GalQ\to\Zell^*$$
 is the $\ell$-adic cyclotomic character, and where $\ell$ is the
 prime of~$\Q$ lying below~$\lambda$
 \cite{\MYTHESIS, Lemma~4.5.1}.  This formula must be modified
 slightly in the case where $E$ is allowed to be a CM-field:

 \proclaim{(3.1) Lemma}
 There is a character of finite order $\epsilon\:\GalQ\to E^*$
 such that $\delta_\lambda= \epsilon\chi_\ell$ for each finite
 prime $\lambda$ of~$E$.  This character is unramified at
 each prime which is a prime of good reduction for~$A$.
 \endproclaim\startproof
 Since the abelian
 representations $\delta_\lambda$
 arise from an abelian variety, they
 have the Hodge-Tate property.  It follows that they
 are locally algebraic
 in the sense of~\cite{\MG, Ch.~III};
 see \cite{\MG, p.~III-50} and \cite{\MYTHESIS, Prop.~1.5.3}.
 One deduces that the family $\delta_\lambda$ is
 associated with an $E$-valued Grossencharacter of type~$A_o$
 of the field~$\Q$ \cite{\MYTHESIS, p.~761}.  Since the
 type-$A_o$ Grossencharacters of~$\Q$ are just products of Dirichlet
 characters with powers of the ``norm" character,
 there is an integer $n$ and an $E$-valued Dirichlet
 character $\epsilon$ so that we have
 $\delta_\lambda= \epsilon\chi_\ell^n$ for each $\lambda$.
 Here, $\ell$ again denotes the residue characteristic of~$\lambda$.
 (We will use the association $\ell\leftrightarrow\lambda$ from
 time to time without comment.)
 We have blurred the distinction
 between characters of $\GalQ$ of
 finite order and Dirichlet characters: if $\epsilon$ is
 a Dirichlet character, then its Galois-theoretic avatar
 takes the value $\epsilon(p)$ on a Frobenius element for~$p$
 in~$\GalQ$.

 By the criterion of N\'eron-Ogg-Shafarevich, $\rho_\lambda$ is
 ramified at a prime $p\ne\ell$ if and only if $A$ has bad
 reduction at~$p$.  Thus, $\delta_\lambda$ is unramified at
 a prime $p\ne\ell$ if $A$ has good reduction at~$p$.  In
 other words, if $p$ is a prime of good reduction for~$A$,
 and if $\ell\ne p$, then
 $\epsilon\chi_\ell^n$ is unramified at~$p$.  Since $\chi_\ell$
 is unramified at~$p$, $\epsilon$ is unramified at~$p$.

 It remains to check that $n=1$.  For this, fix a prime number~$\ell$.
 It is well known that the action of~$\GalQ$
 on $\det_{\Qell}(V_\ell)$ is given by the character
 $\chi_\ell^{\dim(A)} = \chi_\ell^{[E\:\Q]}$.
 On the other hand, the
 determinant of the action of~$\GalQ$ on~$V_\ell$ is
 given by the character
 $$\prod_{\lambda \mid\ell} \Norm_{E_\lambda/\Qell}(\delta_\lambda)
 = \Norm_{E/\Q}(\epsilon)\cdot\chi_\ell^{n\cdot [E\:\Q]}.$$
 (Here, $\Norm$ denotes a norm.)
 Since $\chi_\ell$ has infinite order, we deduce that
 $\Norm(\epsilon)=1$ and that $n=1$.  \endproof

 \proclaim{(3.2) Lemma}
 Each character $\delta_\lambda$ is odd in the sense that it
 takes the value $-1$ on complex conjugations
 in~$\GalQ$.
 \endproclaim\noindent
 [Since $\chi_\ell$ is an odd character, the Lemma may be
 reformulated as the statement that $\epsilon(-1)=+1$.]
 For the proof,
 we use the
 comparison isomorphism
 $$V_\lambda \,\approx\, H_1(A(\C),\Q)\otimes_E E_\lambda$$
 resulting from an embedding $\Qbar \hookrightarrow \C$.
 In this view of $V_\lambda$, the complex conjugation
 of $\GalQ$ acts on $V_\lambda$ as $F_\infty\otimes 1$, where
 $F_\infty$ is the standard ``real Frobenius" on~$H_1(A(\C),\Q)$
 (cf.~\cite{\DELCORVALIS, \S0.2}).
 In particular, $\delta_\lambda$ is odd if and only if we
 have $\det F_\infty = -1$, where the determinant is taken
 relative to the $E$-linear action of
 $F_\infty$ on~$H_1(A(\C),\Q)$.

 Since $F_\infty$ is an involution, and $H_1(A(\C),\Q)$ has
 dimension~two, the indicated determinant is~$+1$ if and only
 if $F_\infty$ acts as a scalar ($=\pm1$)
 on~$H_1(A(\C),\Q)$.  To prove that $F_\infty$ does {\it not\/}
 act as a scalar, we recall that $F_\infty\otimes 1$ permutes
 the two subspaces $H_{0,1}$ and $H_{1,0}$
 of~$H_1(A(\C),\Q)\otimes_\Q\C$
 in the Hodge decomposition of~$H_1(A(\C),\Q)\otimes_\Q\C$.
 \endproof

 \proclaim{(3.3) Proposition}
 For each $\lambda$, $\rho_\lambda$ is an absolutely irreducible
 two dimensional representation of~$\GalQ$ over~$E_\lambda$.
 We have $\End_{\Q_\ell[\GalQ]} V_\lambda = E_\lambda$.
 \endproclaim\startproof
 Faltings \cite\FALT\ proved that
 $V_\ell$ is a
 semisimple $\Qell[\GalQ]$-module whose commutant is
 $E\otimes{\Qell}$.  (This is the Tate conjecture for endomorphisms
 of~$A$.)  Since $V_\ell$ is the product of the~$V_\lambda$ with
 $\lambda\mid\ell$, each module
 $V_\lambda$ is semisimple over $\Qell[\GalQ]$
 and satisfies $\End_{\Qell[\GalQ]}V_\lambda = E_\lambda$.
 This implies that $V_\lambda$ is simple over~$E_\lambda$, and
 that $\End_{E_\lambda} V_\lambda = E_\lambda$.
 The absolute irreducibility follows.
 \endproof

 For each prime $p$ at which $A$ has good reduction, let
 $a_p$ be the element of~$E$ such that
 $$ a_p = \trace_{E_\lambda}(\Frob_p \mid V_\lambda) $$
 whenever $\ell\ne p$.  Let $\overline{\phantom a}$ denote the canonical
 involution on~$E$: this involution is the identity if $E$
 is totally real, and the ``complex conjugation"
 on~$E$ if $E$ is a CM~field.  The involution  $\overline{\phantom a}$
 is the Rosati involution on~$E$ induced by every polarization of
 $A/\Q$.

 \proclaim{(3.4) Proposition}  We have $a_p = {\overline{a}}_p \epsilon(p)$
 for each prime $p$ of good reduction.\endproclaim
 \startproof
 Let $\ell$ be a prime number. Let $\sigma\:E\hookrightarrow\Qellbar$
 be an embedding of fields, and let $\sigmabar$ be the conjugate
 embedding $x\mapsto \sigma(\overline{x})$.  Let $V_\sigma
 = V_\ell \otimes_{E\otimes\Qell}\Qellbar$, where the tensor
 product is taken relative to the map of $\Qell$-algebras
 $E\otimes\Qell\to\Qellbar$ induced by~$\sigma$.
 Define $V_\sigmabar$ similarly, using $\sigmabar$.

 Fix a polarization of~$A/\Q$.  The associated
 $e_{\ell^\nu}$-pairings of Weil induce a bilinear map
 $\langle \, , \, \rangle\:V_\ell \times V_\ell \to \Qell(1)$,
 where $\Qell(1)$ is (as usual) the $\Qell$-adic Tate
 module attached to the $\ell$-power roots of unity in~$\Qbar$.
 We have $\langle ex, y\rangle = \langle x, \overline e x\rangle$
 for $e\in E$ and $x,y\in V_\ell$.
 Further, this pairing is $\GalQ$-equivariant in the sense that
 we have
 $\langle \ggg x, \ggg y\rangle = \ggg\langle x,y\rangle$
 for $g\in\GalQ$.
 After extending scalars from
 $\Qell$ to $\Qellbar$, we find an isomorphism of
 $\Qellbar[\GalQ]$-modules
 $$ V_\sigmabar \approx \Hom(V_{\sigma},\Qellbar(1)).$$
 Now (3.1) implies that
 the determinant of $V_{\sigma}$ is a
 one-dimensional $\Qellbar$-vector space
 on which $\GalQ$ acts by the character $\sigma(\epsilon\chi_\ell) =
 \ssigma\epsilon\chi_\ell$.
 Since $V_{\sigma}$ is of dimension two, this gives
 $$\Hom(V_{\sigma},\Qellbar(\ssigma\epsilon\chi_\ell)) \approx
    V_{\sigma}.$$
 In view of the fact
 that $\Hom(V_{\sigma},\Qellbar(\ssigma\epsilon\chi_\ell))$
 is the twist by~$\ssigma\epsilon$ of
 $\Hom(V_{\sigma},\Qellbar(1))$,
 we get $V_{\sigma}\approx V_\sigmabar(\ssigma\epsilon)$.
 For $p\ne\ell$ a prime of good reduction, the
 trace of $\Frob_p$ acting on $V_\sigma$ is $\sigma(a_p)$,
 and similarly the trace of $\Frob_p$ acting
 on~$V_\sigmabar(\ssigma\epsilon)$ is $\ssigma\epsilon(p)\sigmabar(a_p)
 = \ssigma\epsilon(p)\sigma(\overline{a}_p)$.
 This gives $a_p = \epsilon(p)\overline{a}_p$, as required.\endproof

 \proclaim{(3.5) Proposition}
 Let $S$ be a finite set of prime numbers including the set of primes
 at which $A$ has bad reduction.  Then the
 field $E$ is generated over~$\Q$ by the $a_p$ with $p\not\in S$.
 \endproclaim
 \startproof
 The Proposition follows from
 the Tate Conjecture for endomorphisms of~$A$
 by a simple argument, which
 we now sketch (see \cite{\MYTHESIS, pp.~788--789} for more
 details).
 Choose a prime number $\ell$, and observe that $\End_{\GalQ}V_\ell
 = E\otimes\Qell$ because of Faltings's
 results quoted in the proof of~(3.3).  Let $\overline{V}_\ell =
 V_\ell\otimes\Qellbar$; then as a consequence we have
 $\End_{\Qellbar[\GalQ]}\overline{V}_\ell =E\otimes\Qellbar$.
 In addition, the semisimplicity of the $V_\ell$ implies that
 $\overline{V}_\ell$ is semisimple as a $\Qellbar[\GalQ]$-module.
 For each $\sigma\:E\hookrightarrow \Qellbar$, let $V_\sigma$
 be as above.  Then the $V_\sigma$ are simple $\Qellbar[\GalQ]$-modules,
 and they are pairwise non-isomorphic.  Indeed, they are
 a priori semisimple, but the commutant of their product
 is $\prod_\sigma \Qellbar$.  It follows that their traces
 are pairwise distinct.  Since the trace of $\Frob_p$ acting
 on $V_\sigma$ is $\sigma(a_p)$ (for $p\not\in S\cup \{\ell\}$),
 the Cebotarev Density Theorem implies that the functions
 $p\mapsto \sigma(a_p), p\not\in S\cup \{\ell\}$ are pairwise
 distinct.
 \endproof

 For the next result, fix a finite set $S$ as in~(3.5), and let
 $F$ be the subfield of~$E$ generated by the numbers $a_p^2/\epsilon(p)$
 with $p\not\in S$.

 \proclaim{(3.6) Proposition} The field $F$ is totally real.
 The extension $E/F$ is abelian.\endproclaim\startproof
 Let $\overline{\phantom a}$ again be the canonical complex
 conjugation of~$E$.  For $p\not\in S$, we have
 $${\overline{a}_p^2\over\overline{\epsilon}(p)} =
    {a_p^2\over \epsilon(p)^2}\epsilon(p) $$
 by~(3.4).  The first assertion of the Proposition then
 follows.  For the second, let $t_p=a_p^2/\epsilon(p)$.
 It is clear that $E$ is contained in the extension of~$F$
 (in an algebraic closure of~$E$)
 obtained by adjoining to~$F$ the square roots of all~$t_p$,
 and all roots of unity.  This gives the second assertion.\endproof

 We now consider the reductions of the $\lambda$-adic
 representations~$\rho_\lambda$.  To do this directly,
 replace $A$ by an abelian variety which is
 $\Q$-isogenous to~$A$
 and which has the property that its ring of
 $\Q$-endomorphisms is the ring of integers $\Cal O$
 of~$E$.
 (This process does not change isomorphism classes of
 the $\lambda$-adic representations $\rho_\lambda$.)
 Write simply $A$ for this new abelian variety.
 For each $\lambda$, consider the kernel $A[\lambda]$:
 this is the group of $\Q$-valued points of~$A$ which are
 killed by all elements of the maximal ideal $\lambda$ of~$\Cal O$.
 The action of~$\Cal O$ on~$A[\lambda]$ makes $A[\lambda]$
 into a two-dimensional vector space over the residue field
 $\F_\lambda$ of~$\lambda$.  The $\F_\lambda$-linear action
 of~$\GalQ$ on~$A[\lambda]$ defines a reduction
 $\rholambdabar$ of~$\rho_\lambda$ (cf.~\cite{\MYTHESIS, II.2}).
 (There should be no confusion with the involution $\overline{\phantom a}$
 which appears above.)

 \proclaim{(3.7) Lemma}
 For all but finitely many $\lambda$, the representation
 $\rholambdabar$ is absolutely
 irreducible.
 \endproclaim\noindent
 A result of Faltings \cite{\FW, Theorem 1, page 204}
 implies that
 the following holds for almost all~$\lambda$:
 $A[\lambda]$ is a
 semisimple $\F_\ell[\GalQ]$ module whose commuting
 algebra is $\F_\lambda$.
 The Lemma follows directly from this statement.\endproof
 \head 4. Conjectural connection with modular forms.\endhead\noindent

 We continue the discussion of~\S3, focusing
 on the
 possibility of linking $A$ with modular forms, at
 least conjecturally.  According to~(3.7), $\rholambdabar$
 is absolutely irreducible for almost all~$\lambda$.
 For each $\lambda$ such that
 $\rholambdabar$
 is absolutely irreducible, conjectures
 of Serre \cite{\DUKE, (3.2.3$_?$--3.2.4$_?$)} state that
 $\rholambdabar$ is ``modular" in the sense that it
 arises from
 the space of mod~$\ell$ cusp forms of a specific
 level $N_\lambda$, weight $k_\lambda$, and character
 $\epsilon_\lambda$.

 These invariants
 are essentially constant as functions of~$\lambda$.  Rather than
 study them for all $\lambda$, we will
 restrict attention to those maximal ideals $\lambda$
 which are odd,
 prime to the conductor of~$A$,
 are unramified in~$E$,
 and have degree~one.
 Let $\Lambda$ be the set
 of such ideals $\lambda$ with the property that
 $\rholambdabar$ is absolutely irreducible.
 Lemma 3.7 implies that $\Lambda$ is an infinite set.

 \proclaim{(4.1) Lemma} The levels $N_\lambda$ are bounded as
 $\lambda$ varies in~$\Lambda$.
 \endproclaim\startproof
 It follows from the definition of~$N_\lambda$ that this level
 divides the conductor of the two-dimensional
 $\ell$-adic representation $\rho_\lambda$.
 (To compare the two conductors, one can use the Hilbert Formula
 of~\cite{\OGG, \S I}.)
 The conductor of~$\rho_\lambda$
 divides the conductor of the full $\ell$-adic
 representation $V_\ell(A)$.
 According to results of A.~Grothendieck \cite{\SGA, Cor.~4.6},
 this latter conductor is independent of~$\ell$.  (It is
 by definition the conductor of~$A$.)\endproof

 \proclaim{(4.2) Lemma} For all $\lambda\in\Lambda$, we
 have $k_\lambda=2$.\endproclaim
 \startproof
 Take $\lambda\in\Lambda$.
 The determinant of $\rholambdabar$ is the reduction
 mod~$\lambda$ of $\delta_\lambda$.  By~(3.1), this
 determinant is the product of
 the mod~$\ell$ cyclotomic
 character $\chiellbar$
 and the reduction mod~$\lambda$ of~$\epsilon$.
 The definition of~$\Lambda$ shows that $\ell$ is a prime
 of good reduction for~$A$,
 so that $\epsilon$ is unramified at~$\ell$ (Lemma~3.1).
 Hence, if $I$ is an inertia group for~$\ell$ in~$\GalQ$,
 we have $\det\rholambdabar\mid I = \chiellbar$.

 Further, suppose that $D\subset\GalQ$ is a decomposition
 group for~$\ell$.  It is clear that $\rholambdabar\mid D$
 is finite at~$\ell$ \cite{\DUKE, p.~189}, since $A$ has good
 reduction at~$\ell$.  Indeed,
 the kernel of multiplication by~$\ell$
 on~$A_{\Qell}$ extends to a finite flat group scheme $\Cal G$
 over~$\Zell$
 because of this good reduction \cite{\SGA, Cor.~2.2.9}.
 The Zariski closure of~$A[\lambda]$
 in~$\Cal G$ then prolongs $\rholambdabar$ to a group scheme of
 type~$(\ell,\ell)$ over~$\Zell$.

 By Proposition~4 of~\cite{\DUKE, \S2.8}, we find that $k_\lambda=2$.%
 \endproof

 \proclaim{(4.3) Lemma} For all but finitely
 many $\lambda\in\Lambda$, we have $\epsilon_\lambda=\epsilon$.
 \endproclaim
 \noindent
 One checks easily that $\epsilon_\lambda=\epsilon$ whenever~$\ell$
 is prime to the order of~$\epsilon$.  We omit the details,
 since Lemma~4.3 will not be used below.\endproof

 \proclaim{(4.4) Theorem}
 Let $A$ be an abelian variety over~$\Q$ of
 \GLTWO-type.
 Assume Serre's
 conjecture \cite{\DUKE, (3.2.4$_?$)}
 on representations of~$\GalQ$.
 Then $A$ is isogenous to
 a $\Q$-simple factor
 of~$J_1(N)$, for some $N\ge1$.
 \endproclaim
 \startproof
 Applying \cite{\DUKE, (3.2.4$_?$)} to
 the representations
 $\rholambdabar$ with $\lambda\in\Lambda$, we find
 that each $\rholambdabar$
 arises from a newform of weight $k_\lambda=2$ and
 level dividing $N_\lambda$.
 Since the $N_\lambda$'s
 are bounded (Lemma~4.1),
  there are only a finite number of such
 newforms.

 Hence there is a fixed newform $f=\sum a_nq^n$ which gives
 rise to an infinite number of the $\rholambdabar$'s.
 Explicitly, we have the following situation.  Let $R$ be
 the ring of integers of the field $\Q(\ldots, a_n, \ldots)$.
 For an infinite number of~$\lambda\in\Lambda$, there is
 a ring homomorphism $\varphi_\lambda\: R \to \Flambdabar$
 mapping $a_p$ to $\trace(\rholambdabar(\Frob_p))$
 for all but finitely many primes~$p$.

 Let $N$ be the level of~$f$, and let
 $A_f$ be the quotient of~$J_1(N)$ which is
 associated to~$f$.
 Let $\lambda$ be a prime for which there is a
 $\varphi_\lambda$ as above.
 By the Cebotarev Density Theorem, we have
 $$A_f[\ell]\otimes_{R/\ell R}\Flambdabar
      \,\approx\, A[\lambda]\otimes_{\F_\lambda}\Flambdabar,$$
 where $\Flambdabar$ is regarded as an $R/\ell R$-module
 via~$\varphi_\lambda$.
 Since $A_f[\ell]\otimes_{R/\ell R}\Flambdabar$ is a quotient
 of $A_f[\ell]\otimes_\Fell\Flambdabar$, and since $\F_\lambda=\Fell$,
 we have $$\Hom_{\Flambdabar[\GalQ]}(A_f[\ell]\otimes_\Fell\Flambdabar,
   A[\lambda]\otimes_\Fell\Flambdabar) \ne 0.$$
 It follows that $\Hom_{\Fell[\GalQ]}(A_f[\ell], A[\lambda])\ne0$.

 Hence we have $\Hom_{\Fell[\GalQ]}(A_f[\ell], A[\ell])\ne0$ for
 an infinite number of prime numbers~$\ell$.
 By the theorem of Faltings quoted above, we get
 $\Hom_\Q(A_f,A)\ne0$.  Since $A$ is a simple abelian variety over~$\Q$,
 $A$ must be a quotient of~$A_f$.  Since $A_f$ is, in turn,
 a quotient of~$J_1(N)$, we deduce that $A$ is a quotient of~$J_1(N)$.
 \endproof

 \head 5. Decomposition over~$\Qbar$.\endhead\noindent
 Suppose again that $A/\Q$ is an abelian variety
 of \GLTWO-type, and let $$E=\Q\otimes\End_\Q(A).$$
 Let $\Cal X = \Q\otimes\End_{\Qbar}(A)$ be the algebra of
 {\it all\/} endomorphisms of~$A$.

 \proclaim{(5.1) Proposition {\rm [Shimura]}}
 Suppose that $A_{\Qbar}$ has a non-zero abelian
 subvariety of CM-type.
 Then $A_{\Qbar}$ is isogenous to a power of a CM elliptic
 curve.
 \endproclaim
 \noindent
 This is Proposition~1.5 in \cite\CFHO.
 For a generalization of this result, see
 Proposition~5.2 of~\cite\HIDA.
 See also the
 discussion in~\S4 of~\cite\BONN.

 \proclaim{(5.2) Proposition}
 Suppose that $A_{\Qbar}$ has no non-zero abelian
 subvariety of CM type.
 Then the center of~$\Cal X$ is a subfield $F$ of~$E$.
 The
 algebra $\Cal X$ is isomorphic either to a matrix ring over~$F$,
 or else to a ring of matrices over a % totally indefinite???
 quaternion division algebra over~$F$.
 \endproclaim
 \startproof
 We employ the same arguments used to prove~(2.1) above
 and Theorem~2.3 of \cite\ANNALS.  Let $D$ be
 the commutant of~$E$ in~$\Cal X$.  Clearly, $D$ is a division
 algebra: otherwise, we can make $E$ act on a proper
 non-zero
 abelian subvariety of $A_{\Qbar}$, contrary to the
 hypothesis that no abelian subvariety of~$A_{\Qbar}$ is
 of CM~type. Also, $E$ is a subfield of~$D$; it is a maximal
 commutative subfield because $A_{\Qbar}$ does not have
 complex multiplication.  Hence $E$ is its own commutant in~$D$.
 Since $D$ is the commutant of~$E$, we get $D=E$.  In particular,
 the center of~$\Cal X$ is contained in~$E$, so that the
 center is a subfield $F$ of~$E$.

 Since $\Cal X$ is now a central simple algebra over~$F$, we have
 $\Cal X \approx M(n,D)$, where $D$ is a division algebra with
 center~$F$, and $n$ is a positive integer.  Suppose that $D$ has
 dimension $t^2$ over~$F$; then $[E\:F]=nt$ since $E$ is a maximal
 commutative subalgebra of~$\Cal X$.  Up to isogeny, $A_{\Qbar}$
 is of the form $B^n$, where the endomorphism algebra of~$B$
 contains $D$.  Following an idea of J.~Tunnell
 (cf.~\cite{\DPP, Th.~1}), we
 fix an embedding $\Qbar\hookrightarrow\C$ and
 form the cohomology group $H^1(B(\C),\Q)$.  This is a $\Q$-vector
 space of dimension
 $$ 2\dim(B) = {2\over n}\dim(A) = {2nt\over n}[F\:\Q]$$
 with a functorial action of~$D$.  Hence the
 $\Q$-dimension $t^2[F\:\Q]$ of~$D$ divides the dimension over~$\Q$
 of~$H^1(B(\C),\Q)$, which is
 $2t[F\:\Q]$. Thus $t\le 2$, so that either  $D=F$, or else
 $D$ is a quaternion division algebra with center~$F$.
 \endproof

 As in Proposition~3.6,
 let $S$ be a finite set of primes containing the primes of
 bed reduction for~$A$.  Let $F$ be the center of~$\Cal X$.
 Then we have:

 \proclaim{(5.3) Theorem}
 The field $F$ is generated by the numbers $a_p^2/\epsilon(p)$
 with $p\not\in S$.
 \endproclaim
 \noindent In other words, the field $F$ which appears in (3.6)
 is the same as the field $F$ in~(5.2).

 Let $\ell$ be a prime which
 splits completely in~$E$, so that all
 embeddings $E\hookrightarrow\Qellbar$ take values in~$\Qell$.
 Choose a finite extension $K$ of~$\Q$ such that all endomorphisms
 of~$A$ are defined over~$K$, and let $H$ be the
 corresponding open subgroup of~$\GalQ$.
 Replacing $H$ be a smaller subgroup of~$\GalQ$ if necessary,
 we may assume that $H$ is contained in the kernel of~$\epsilon$.
 We have $\Cal X\otimes\Qell = \End_{\Qell[H]}{V}_\ell$, by
 Faltings's results \cite{\FALT}.
 The center of~$\Cal X\otimes\Qell$ is $F\otimes\Qell$.

 The Tate module ${V}_\ell$ decomposes as
 a product $\prod_\sigma V_\sigma$, where
 $\sigma$ runs over the set $\Sigma$ of
 embeddings $\sigma\:E\to \Qell$,
 and where $V_\sigma=V_\ell\otimes_{E\otimes\Qell}\Qell$,
 with $\Qell$ being regarded as an $E\otimes\Qell$ module
 via~$\sigma$.
 (Cf.~the proof of~(3.4).)
 Each
 $V_\sigma$ is a simple $\Qell[H]$-module because $A$ has
 no CM subvariety and because the action of~$H$ on
 $V_\ell$ is semisimple (Faltings).  Hence $\End_H V_\sigma =\Qell$
 for each~$\sigma$.

 For each prime $v$ of~$K$ which is prime to $\ell$ and the set
 of bad primes for~$A$, there is a ``trace of Frobenius" $t_v\in E$
 associated with~$v$.  We have $\trace(\Frob_v\mid V_\sigma) =
 \sigma(t_v)$ for each~$\sigma$.  One finds for $\sigma,\tau\in\Sigma$
 that $V_\sigma$ and $V_\tau$ are isomorphic $\Qell[H]$-modules
 if and only if $\sigma(t_v)=\tau(t_v)$ for all~$v$, i.e., if
 and only if $\sigma|_L = \tau|_L$, where $L=\Q(\ldots,t_v,\ldots)$
 (cf.~\cite{\MYTHESIS, \S IV.4}).
 This implies that the center of~$\Cal X\otimes\Qell$ is~$L\otimes\Qell$.
 Thus $F\otimes\Qell=L\otimes\Qell$
 (equality inside $E\otimes\Qell$), which implies that~$F=L$.

 Suppose now that $\sigma=\tau$ on~$L$, so that $V_\sigma\approx
 V_\tau$ as representations of~$H$.
 A simple argument shows that
 there is a
 character $\phi\:\GalQ\to\Q^*_\ell$ such that $V_\sigma \approx
 V_\tau\otimes\phi$ as $\Qell[\GalQ]$-modules.
 The character $\phi$
 is necessarily
 unramified at all primes $p\ne\ell$ which are
 primes of
 good reduction for~$A$.
 Taking traces, we get
 $\sigma(a_p) = \phi(p)\tau(a_p)$ for all such primes.
 By considering determinants, we get
 $\ssigma\epsilon= \phi^2\cdot\ttau\epsilon$.
 These two equations show that $\sigma$ and~$\tau$ agree on
 $a_p^2/\epsilon(p)$ for all good primes $p\ne\ell$.

 It follows by Galois theory that we have $a_p^2/\epsilon(p)\in L$
 for all good reduction primes $p\ne\ell$.  Hence $a_p^2/\epsilon(p)\in F$
 for all such~$p$.  By varying~$\ell$ we obtain the inclusion
 $\Q(\ldots,a_p^2/\epsilon(p),\ldots) \subseteq F$, where $p$ runs
 over the set of all primes of good reduction for~$A$.

    To prove the opposite inclusion, we must show that if
 $\sigma(a_p^2/\epsilon(p)) = \tau(a_p^2/\epsilon(p))$ for all~$p$,
 then $V_\sigma$ and $V_\tau$ are $H$-isomorphic.
 By the Cebotarev Density Theorem, the hypothesis implies that
 the functions on~$\GalQ$ ``$\trace^2/\det$" are the same for
 $V_\sigma$ and~$V_\tau$.  In particular, we have
 $\trace(h| V_\sigma) = \pm \trace(h | V_\tau)$ for
 all $h\in H$.  This equality implies that $V_\sigma$ and $V_\tau$
 become isomorphic after $H$ is replaced by an open subgroup $H_o$
 of~$H$ (cf.~\cite{\SERREINV, p.~324}).  Since $H$ was already
 chosen ``sufficiently small," we find that $V_\sigma$
 and $V_\tau$ are indeed isomorphic as $\Qell[H]$-modules.\endproof

 \proclaim{(5.4) Corollary} The center $F$ of~$\Cal X$ is a
 totally real number field.  The extension $E/F$ is abelian.
 \endproclaim\startproof
 The Corollary follows from (5.3) and~(3.6).
 \endproof

 Let $g$ be an element of~$\GalQ$, and consider the automorphism
 $x\mapsto \ggg x$
 of~$\Cal X$ which is induced by~$g$.
 This automorphism is necessarily inner (Skolem-Noether theorem),
 and it fixes $E$; therefore it is given by conjugation by an
 element $\alpha(g)$
 of~$E^*$ which is well defined modulo $F^*$.  The map
 $\alpha\:\GalQ\to E^*/F^*$ is a continuous homomorphism
 (cf.~\cite{\DPP, p.~268}).  It is a fact that $\alpha$ is
 unramified at all primes $p$ at which $A$ has semistable reduction
 \cite{\ANNALS, Th.~1.1}.

 \proclaim{(5.5) Theorem}
 We have $\alpha^2\equiv\epsilon$ {\rm (mod~$F^*$)}. Moreover,
 suppose that $p$ is a prime of good reduction for~$A$
 such that $a_p\ne0$.  Then $\alpha(\Frob_p)\equiv a_p$
 {\rm (mod~$F^*$)}.
 \endproclaim
 \startproof To prove
 the first assertion, let $\ell$ be a prime which
 splits completely in~$E$.  We must prove that
 $\sigma(\alpha^2(g)/\epsilon(g))=
 \tau(\alpha^2(g)/\epsilon(g))$ whenever $\sigma$ and~$\tau$ are embeddings
 $E\hookrightarrow \Qell$ which agree on~$F$.

 If $\sigma$ and~$\tau$ have this property, then there is
 a character $\phi\:\GalQ\to\Q^*_\ell$ for which $V_\sigma\approx
 V_\tau\otimes\phi$ (the notation is as in the proof of~(5.3)).
 The Galois group $\GalQ$ acts on $\Hom(V_\sigma,V_\tau)$
 by multiplication by~$\phi(g)^{-1}$, but also by conjugation by~$\alpha(g)$.
 Since $\alpha(g)$ acts on~$V_\sigma$ and $V_\tau$ by
 $\ssigma\alpha(g)$ and $\ttau\alpha(g)$ (respectively),
 $\alpha(g)$ acts on $\Hom(V_\sigma,V_\tau)$ by
 $\ttau\alpha(g)/\ssigma\alpha(g)$.
 Hence, $\phi(g)=\ssigma\alpha(g)/\ttau\alpha(g)$
 in~$\Qell$.
 On the other hand, as remarked during the proof of~(5.3),
 we have $\ssigma\epsilon=\phi^2\cdot\ttau\epsilon$.
 The two equalities give the required conclusion
 $$ \sigma(\alpha^2/\epsilon) = \tau(\alpha^2/\epsilon).$$

 For the second assertion, we choose $\ell\ne p$.
 With the notation as above,
 $\phi(p) = \ssigma\alpha(\Frob_p)/\ttau\alpha(\Frob_p)$.
 However, we have $\sigma(a_p) = \phi(p)\tau(a_p)$, as
 noted during the proof of~(5.3). On comparing two
 formulas for~$\phi(g)$, we get $\sigma(\alpha(\Frob_p)/a_p)=
 \tau(\alpha(\Frob_p)/a_p)$.\endproof

 \proclaim{Remark}\rm
 It should be easy to prove
 that the set of~$p$ for which $a_p=0$ has
 density~0, since the image of~$\GalQ$ in $\Aut V_\sigma$
 is open in $\Aut V_\sigma$, for any~$\sigma$.
 \endproclaim
 \noindent
 Let $\tilde\alpha$ be a lift of~$\alpha$ to a function
 $\GalQ\to E^*$.
 The function
 $$c\:(g_1,g_2)\mapsto {\tilde\alpha(g_1)\tilde\alpha(g_2)
  \over\tilde\alpha(g_1g_2)}$$
 is then a 2-cocycle on~$\GalQ$ with values in~$F^*$.  (We regard
 $F^*$ as a $\GalQ$-module with trivial action.)
 The image $[c]$ of~$c$ in $H^2(\GalQ,F^*)$ is independent
 of the choice of~$\tilde\alpha$.  Since $\alpha^2$ lifts to the
 character $\epsilon$, it is clear that $[c]$
 has order dividing~two.

 Consider the map
 $$ \xi\:H^2(\GalQ,F^*) \, \longrightarrow \, H^2(\Gal(\Qbar/F),\Qbar^*)
            = \Br F$$
 obtained from the restriction
 map $H^2(\GalQ,F^*)\to H^2(\Gal(\Qbar/F),F^*)$
 and the map $H^2(\Gal(\Qbar/F),F^*) \to H^2(\Gal(\Qbar/F),\Qbar^*)$
 induced by the inclusion of $F^*$ into $\Qbar^*$.
 (We view $\Qbar^*$ as a $\Gal(\Qbar/F)$-module in the
 standard way; i.e., we use the Galois action.  The notation
 ``$\Br F$" indicates the Brauer group of~$F$.)

 \proclaim{(5.6) Theorem {\rm [Chi]}}
 The class of~$\Cal X$ in~$\Br F$ coincides with $\xi([c])$.
 \endproclaim
 \startproof
 According to \cite{\CHI, Th.~3.4}, the algebra $\Cal X$
 is isomorphic to a twisted matrix algebra $(\End_F E)(\alpha)$.
 Theorem~4.8 of~\cite\CHI\ expresses $[(\End_F E)(\alpha)]\in\Br F$
 as the image of a certain two-cocycle whose values are Jacobi
 sums.
 By \cite{\DPP, Prop.~1}, the class of this Jacobi-sum
 cocycle coincides with the class of the two-cocycle
 $$ (g,h)\mapsto \tilde\alpha(h)^{g-1}
   = {\ggg\tilde\alpha(h)\over\tilde\alpha(h)}$$
 on $\Gal(\Qbar/F)$.  We wish to compare this with the two-cocycle
 $$ (g,h)\mapsto
      {\tilde\alpha(g)\tilde\alpha(h)\over\tilde\alpha(gh)}.$$
 The product of the two is the map
 $$ (g,h)\mapsto
     {\ggg\tilde\alpha(h)\tilde\alpha(g)\over\tilde\alpha(gh)},$$
 which is a coboundary.  Hence the class of~$\Cal X$ in~$\Br F$
 is the negative of $\xi([c])$.  Since $[c]$ (and
 $[\Cal X]$) have order two, we
 get the Theorem as stated.\endproof

 A final remark about $\Cal X$ concerns the isogeny class
 of the simple factors of~$A_{\Qbar}$.  We have $\Cal X\approx
 M(n,D)$ for some positive integer~$n$, where $D$ is a
 division algebra of dimension one or~four over~$F$.
 Accordingly, $A_{\Qbar}$ decomposes up to isogeny as
 the product of $n$ copies a
 simple abelian variety $B$ over~$\Qbar$
 whose endomorphism algebra is isomorphic to~$D$.
 Since $A$ is defined over~$\Q$, we have
 $$ \ggg B\times\cdots\times\ggg B \,\,\sim \,\,\ggg A_{\Qbar}
  \,\,\sim A_{\Qbar} \,\,\sim B\times\cdots\times B$$ for each $g\in\GalQ$.
 (The sign ``$\sim$" indicates an isogeny.)
 By the uniqueness of decomposition up to isogeny, we have
 $\ggg B\sim B$.  One says that the isogeny class of~$B$
 is ``defined over~$\Q$."

 Consider the special case where $n=\dim A$, so that $B$ is
 of dimension~one.  The elliptic curve $B$ is then
 isogenous to
 each of its conjugates (over~$\Qbar$).  Borrowing
 (and bending) a
 term used by B.~H.~Gross \cite\BHG, we say that
 $B$ is a ``$\Q$-curve."

 \head 6. $\Q$-curves as factors of abelian varieties of \GLTWO-type.%
 \endhead\noindent
 Suppose that $C$ is an elliptic curve over~$\Qbar$ which occurs
 as a simple factor of a primitive abelian variety of \GLTWO-type
 over~$\Q$ with no complex multiplication over~$\Qbar$.
 As we have just seen, $C$ is a $\Q$-curve; it is, more
 precisely, a $\Q$-curve with no complex multiplication.

 In this \S, we prove the converse:

 \proclaim{(6.1) Theorem}
 Suppose that $C$ is an elliptic curve over~$\Qbar$ with no
 complex multiplication.  Assume that $C$ is isogenous to each
 of its conjugates $\ggg C$ with $g\in\GalQ$.
 Then there is a primitive abelian variety $A$
 of \GLTWO-type over~$\Q$ such that $C$ is a simple factor
 of~$A$ over~$\Qbar$.
 \endproclaim

 \proclaim{(6.2) Corollary}
 Suppose that $C$ is as in~(6.1).
 Assume Serre's
 conjecture \cite{\DUKE, (3.2.4$_?$)}
 on representations of~$\GalQ$.
 Then $C$ is a simple factor, over~$\Qbar$, of $J_1(N)$ for
 some $N\ge1$.
 \endproclaim
 \noindent
 The Corollary follows directly from (6.1) and~(4.4).\endproof
 \medskip
 \noindent
 {\sl Proof of~(6.1).}\enspace
 Let $C$ be as in the statement of the Theorem.  We can find
 a model $C_o$
 of~$C$ over a number field $K\subset\Qbar$.
 We may assume that $K/\Q$ is a Galois extension.
 For each $g\in\Gal(K/\Q)$, $C_o$ and $\ggg C_o$ are
 $\Qbar$-isogenous.  Enlarging $K$ if necessary, we may assume
 that there are isogenies $\mu_g\: \ggg C_o\to C_o$ defined
 over~$K$.

 The map
 $$ c\:(g,h)\mapsto \mu_g \ggg \mu_h \mu^{-1}_{gh}$$
 may be regarded as $\Q^*$-valued, since $\Q\otimes\End_K(C_o) =\Q$.
 A short
 computation shows that $c$ is a two-cocycle on~$\Gal(K/\Q)$
 with values in~$\Q^*$.  The class of~$c$ in $H^2(\Gal(K/\Q),\Q^*)$
 is independent of the choices of the $\mu_g$.
 By inflation, we may (and will) regard $c$ as a locally
 constant function on~$\GalQ$.

 \proclaim{(6.3) Theorem {\rm [Tate]}}
 Let $M$ be the discrete $\GalQ$-module $\Qbar^*$ with
 trivial action.  Then $H^2(\GalQ,M)=0$.
 \endproclaim\noindent
 This theorem of Tate is proved as Theorem~4
 in Serre's article \cite\JPSDURHAM,
 with $\Qbar$ replaced by~$\C$.
 The proof exposed by~Serre works equally well in the
 present case.\endproof

 Because of Tate's theorem, there is a locally constant
 function $\alpha\:\GalQ\to\Qbar^*$ such that
 we have the identity among functions $G\times G\to\Qbar^*$
 $$ c(g,h)= {\alpha(g)\alpha(h)\over\alpha(gh)}.$$
 After again enlarging $K$, we may identify $\alpha$
 with a function on~$\Gal(K/\Q)$ and regard $g$ and~$h$
 as elements of this Galois group.

 Let $E$ be the extension of~$\Q$ generated by the values
 of~$\alpha$.  The definition of~$c$ shows that we have
 $$ c(g,h)^2 = {\deg\mu_g \deg\mu_h\over\deg\mu_{gh}},$$
 where ``deg" denotes the degree of an isogeny between
 elliptic curves (or, more generally, of a non-zero element
 of $\Q\otimes\Hom(C_1,C_2)$, where $C_1$ and $C_2$ are
 elliptic curves).  It follows that the function
 $$ \epsilon\:g\mapsto {\alpha^2(g)\over\deg\mu_g}$$
 is a Dirichlet character $\Gal(K/\Q)\to E^*$.
 Because $\alpha^2 \equiv \epsilon$ (mod~$\Q^*$), the
 field $E$ is an abelian extension of~$\Q$ (cf.~Prop.~3.6).

 If one performs the analysis of~\S3 on the abelian
 variety~$A$ which we construct below, one finds (e.g., in~(3.1))
 another Dirichlet character called $\epsilon$.
 It is very likely that the two~$\epsilon$'s are equal
 up to ``sign," i.e., possible inversion (cf.~Lemma~7.1 below).

    To simplify notation, let us write simply $C$ for the elliptic
 curve $C_o$ over~$K$ and $C_{\Qbar}$ for the curve originally
 called~$C$.  Let $B$ be the abelian variety $\Res_{K/\Q} C$,
 where ``Res" is Weil's ``restriction of scalars"
 functor~\cite{\WEIL, \MILNE}.
 Then $B$ is an abelian variety over~$\Q$ of dimension $[K:\Q]$.
 It represents the functor on~$\Q$-schemes
 $S\mapsto C(S_{K})$; in particular, we have
 $\Hom_\Q(X,B) = \Hom_K(X_{K},C)$ whenever $X$ is an abelian
 variety over~$\Q$.

 Applying this formula in the special case where $X=B$,
 we find $\End_\Q(B)=\Hom_K(B_{K},C)$.
 On the other hand
 \cite{\WEIL, p.~5},
 $$B_{K} = \prod_{\sigma\in\Gal(K/\Q)} \ssigma C.$$
 Hence we have
 $$ \Q\otimes\End_Q(B) = \prod_\sigma \Q\otimes\Hom_K(\ssigma C,C).$$
 Since $\Q\otimes\Hom(\ssigma C,C)$ is the one-dimensional
 vector space generated by $\mu_\sigma$, we find that
 $\Cal R\seteq\Q\otimes\End_\Q(B)$
 may be written as $\prod \Q\cdot\mu_\sigma$.
 Let $\lambda_\sigma$ be the element of $\Cal R$ which
 corresponds to
 $\mu_\sigma\:\ssigma C\to C$; then $\Cal R$ is a $\Q$-algebra
 with vector space basis $\lambda_\sigma$.

 \proclaim{(6.4) Lemma} We have $\lambda_\sigma\lambda_\tau
 = c(\sigma,\tau)\lambda_{\sigma\tau}$ in~$\Cal R$ for
 $\sigma,\tau\in\Gal(K/\Q)$.
 \endproclaim
 \startproof
 Each map $\lambda_\sigma$ acts on $B_{K}=\prod_g \ggg C$ as a
 ``matrix": it sends the factor
 ${}^{g\sigma}\mkern-\thinmuskip C$ to~$\ggg C$
 by $\ggg\mu_\sigma$ (cf.~\cite{\BHG. \S15}).
 A short computation using the identity
 $$ \mu_\sigma\ssigma\mu_\tau = c(\sigma,\tau)\mu_{\sigma\tau}$$
 gives the desired formula.
 \endproof

 The algebra $\Cal R$ is thus a ``twisted group algebra"
 $\Q[\Gal(K/\Q)]$ in which we have the multiplication table
 $[\sigma][\tau]=c(\sigma,\tau)[\sigma\tau]$.
 The obvious map of $\Q$-vector spaces
 $$\omega\:\Cal R\to E, \qquad \lambda_\sigma\mapsto \alpha(\sigma)$$
 is in fact a surjective homomorphism of $\Q$-algebras
 because of~(6.4).

 In \cite\GRBU, $\Q\otimes\End_\Q(B)$ is studied in an analogous
 situation where $C$ has complex multiplication, and where
 $K$ is the
 Hilbert class field of the field of complex multiplication.
 (The curve $C$ is then a ``$\Q$-curve" in the original sense
 of the term.)  Criteria are given for $\Q\otimes\End_\Q(B)$
 to be (i) a commutative algebra, i.e., a
 product of number fields, and (ii) a product of
 {\it totally real\/} fields.  It might be interesting to formulate
 similar criteria in our context.

 Let $T$ be the abelian variety $\prod_\sigma C$ over~$K$.  We
 write $C_\sigma$ for the copy of~$C$ in the $\sigma$th place,
 so that $T$ becomes the product $\prod C_\sigma$.
 We have a map $\Cal R\to \Q\otimes\End T$
 given as follows: For $g\in\Gal(K/\Q)$, the element
 $\lambda_g = [g]$ of~$\Cal R$ acts on~$T$ by sending
 $C_\sigma$ to~$C_{g\sigma}$
 by the map (of elliptic curves up to isogeny)
 ``multiplication by~$c(g,\sigma)$."  One checks directly that
 this is a homomorphism of $\Q$-algebras.  The variety~$T$ will
 be interpreted as $\Cal R\otimes_\Q C$ by readers who are
 fond of such tensor-product constructions.

 Let $\iota\:T\,\buildrel\sim\over\to\,
      B_{K} = \prod \ssigma C$ be the isomorphism
 of abelian varieties up to isogeny which takes the factor
 $C_\sigma$ of~$T$ to the factor ${}^{\sigma^{-1}}\mkern-\thinmuskip
 C$ of~$B_{K}$, via
 the map ${}^{\sigma^{-1}}\mkern-\thinmuskip\mu_\sigma$.

 \proclaim{(6.5) Proposition} The map $\iota$ is $\Cal R$-equivariant,
 for the action of~$\Cal R$ on~$T$ just defined and for the
 structural action of~$\Cal R$ on~$B$.
 \endproclaim\startproof
 The proof of this Proposition is an uninteresting computation,
 which is omitted.

 \proclaim{(6.6) Corollary} The Lie algebra $\Lie(B/\Q)$ is
 a free $\Cal R$-module of rank~one.\endproclaim
 \startproof
 The statement to be proved is true if and only if $\Lie(B_{K}/K)$
 is free of rank~one over $R\otimes_\Q K$.
 By~(6.5), $\Lie(B_{K}/K)$ may be identified with $\Cal R\otimes_\Q
 \Lie(C/K)$, with $\Cal R$ operating trivially on the second factor.
 Hence $\Lie(B_{K}/K)$ is indeed free of rank~one over
 $\Cal R\otimes_\Q K$.
 \endproof

    To complete the proof of~(6.1), we let $A$ be the abelian
 variety $E\otimes_{\Cal R} B$, where $E$ is viewed as a
 $\Cal R$-module via~$\omega$.  Explicitly,
 use the fact that $\Cal R$ is a semisimple
 $\Q$-algebra to write
 $\Cal R$ as a direct sum of its quotient $E$ with
 the kernel of the map $\omega$.  Let $\pi\in\Cal R$ be
 the projector onto~$E$, and let $A\subseteq B$ be the image
 of~$\pi$, viewed as an endomorphism of~$A$ up to isogeny.
 (In other words, $A$ is the image of $m\cdot\pi$, where
 $m$ is a positive integer chosen so that $m\cdot\pi$ is
 a true endomorphism of~$B$.)
 Then $A$ is an abelian subvariety of~$B$, defined over~$\Q$,
 whose algebra of $\Q$-endomorphisms is~$E$.
 By~(6.6), $E$ acts without multiplicity on~$\Lie(B)$ and
 therefore, in particular, without multiplicity on~$\Lie(A)$.
 Hence $A$ has dimensional equal to~$[E\:\Q]$, which means
 that $A$ is of~\GLTWO-type.  (It is clear that $A$ is non-zero
 because $\pi$ is non-zero.)

 Since $B_K$ is isogenous to a product of copies of~$C$,
 the same holds true for~$A_K$.  Thus $C$ is
 a quotient of~$B_K$.
 \endproof

 \head 7. $\Q$-curves over quadratic fields.%
 \endhead\noindent

 Suppose that $C$ is a $\Q$-curve as above and that
 $K$ is a
 quadratic field.  Let $\sigma$ be the non-trivial
 automorphism of $K$ over~$\Q$.
 Then, by hypothesis, there is a $K$-isogeny
 $\mu=\mu_\sigma\:\ssigma C\to C$.
 We take the identity map
 for $\mu_1$, where ``1" is the identity automorphism of $K$.
 The cocycle~$c$ takes the value~1 on all elements of
 $\Gal(K/\Q)\times\Gal(K/\Q)$ other than $(\sigma,\sigma)$.
 Its value on that pair is the non-zero
 integer $m$ such that
 $\mu\circ\ssigma\mu$ is multiplication by $m$ on~$C$.
 The algebra~$\Cal R$ may be written $\Q[X]/(X^2-m)$, where
 $X$ corresponds to the element of~$\Cal R$ we have
 been calling~$[\sigma]$.

 Let us split $c$ by defining $\alpha\:\Gal(K/\Q)\to\Qbar^*$
 to be the map taking $1$ to~$1$ and $\sigma$ to a square
 root of~$m$.  The character
 $$\theta\: g\mapsto {\alpha^2(g)\over\deg\mu_g}$$
 is then trivial if $m$ is positive and of order two
 if $m$ is negative.  In the case where $\theta$ is of order
 two, it is an isomorphism $\Gal(K/\Q)\buildrel\sim\over\to \{ \pm 1\}$.

 If $m$ is a perfect square, then we have $E=\Q$ in the
 notation of~\S6.  The abelian variety $A$ is then
 a model of~$C$ over~$\Q$.

 {\it Assume now that $m$ is not a
 perfect square.}  Then $\Cal R=E$ is a quadratic number field,
 and we have $B=A$ in the notation of~\S6.
 The field $E$ is real if $m$ is positive and imaginary
 if $m$ is negative.  Thus $E$ is real if and only if $\theta$
 is trivial.

 The $\lambda$-adic representations of~$A$
 define a Dirichlet character $\epsilon$ (Lemma~3.1).
 According to~(3.2), $\epsilon$ is an even
 character. Also, $\epsilon$ is non-trivial if and only
 if $E$ is imaginary~(3.4).  Thus $\epsilon$ if non-trivial
 if and only if $\theta$ is non-trivial.

 \proclaim{(7.1) Lemma}
 The characters $\epsilon$ and $\theta$ are equal.\endproclaim\startproof
 By~(6.5), $A_K$ is $K$-isogenous to the abelian variety
 ``$E\otimes C$," i.e., to the product of two copies of~$C$ with
 $E$ acting through a regular representation $E\hookrightarrow M(2,\Q)$.
 In particular, the $\lambda$-adic representations of~$A_K$ are just
 the $\ell$-adic representations of~$C$, viewed as taking values in
 $\bold{GL}(2,E_\lambda)$ rather than in~$\bold{GL}(2,\Qell)$.
 This implies that the determinants of the $\rho_\lambda |_{\Gal(\Qbar/K)}$
 are the cyclotomic characters~$\chi_\ell$.
 Hence $\epsilon$ is trivial on~$\Gal(\Qbar/K)$, and therefore
 may be identified with a character of~$\Gal(K/\Q)$, which
 has order~two.
 Since $\theta$ is also a character of this latter group, and since
 the two characters are simultaneously non-trivial, they are
 equal.\endproof

 As mentioned above,
 it seems very likely that Lemma~7.1
 (quite possibly in the form $\epsilon=\theta^{-1}$)
 generalizes to the situation
 of~\S6.

 \proclaim{(7.2) Proposition {\rm [Serre]}}
 At least one of the two quadratic fields $E$, $K$ is a
 real quadratic field.\endproclaim
 \startproof
 We give two proofs, the first of which was communicated to the
 author by Serre:
 Assume that $K$ is a complex quadratic field.
 After we
 embed $K$
 in~$\C$, the automorphism $\sigma$ of~$K$
 becomes the restriction to~$K$ of $\bar{\phantom L}$,
 complex conjugation on~$\C$.
 Choose a holomorphic differential  $\omega$
 on~$C$.
 Then
 $C_\C$ may be identified with
 the curve $\C/L$, with $L\subset\C$ the period lattice
 of~$\omega$.  The curve $\ssigma C$ becomes the complex
 conjugate $\C/\bar L$
 of~$C$.  The isogeny
 $\mu$ is induced by the map ``multiplication by $\gamma$"
 on~$\C$, for some non-zero complex number~$\gamma$.
 The integer $m$ may be identified with
 $\gamma\bar\gamma$, and is therefore positive.
 This means that $E$ is a real quadratic field.

 Another proof can be given as follows.  Assume that $E$ is
 an imaginary quadratic field.  Then $\epsilon$ is a
 non-trivial character of~$\Gal(K/\Q)$.  This character
 is {\it even\/} by~(3.2).
 Therefore $K$ is imaginary.
 \endproof

 The case where $E$ is imaginary, so that $K$ is real,
 was treated in detail
 by Shimura (see \cite{\CFHO, \S10} for numerical examples).
 This case was later studied by Serre \cite{\DUKE, p.~208},
 who pointed out that Theorem~4.4 holds in this context.

 \head 8. Descent of abelian varieties up to isogeny.%
 \endhead\noindent
 Suppose that $L/K$ is a Galois extension of~fields and
 that $A$ is an abelian variety over~$L$.  A well-known
 theorem of Weil \cite\WEILFD\ states
 that $A$ has a model over $K$ if and only if there are
 isomorphisms $\mu\:\ssigma A\buildrel\sim\over\to A$ ($\sigma\in\Gal(L/K)$)
 which satisfy the compatibility condition
 $$\mu_\sigma\ssigma\mu_\tau = \mu_{\sigma\tau}.\leqno(8.1)$$

 As an application of the techniques encountered in~\S6, we will
 prove an analogous criterion
 for abelian varieties {\it up to isogeny}.
 Namely, suppose that $A$ is isogenous (over~$L$) to an abelian
 variety $B/L$ which has a model over~$K$.  Then one finds
 isomorphisms $\mu_\sigma\:\ssigma A\buildrel\sim\over\to A$
 of abelian varieties up to isogeny over~$L$ which satisfy the
 compatibility condition~(8.1).  Conversely, one has
 \proclaim{(8.2) Theorem}
 Suppose that there are isomorphisms of abelian varieties up
 to isogeny over~$L$,  $\mu\:\ssigma A\buildrel\sim\over\to A$,
 which satisfy~(8.1).  Then there is an abelian variety $B$ over~$K$
 such that $A$ is $L$-isogenous to $B_L$.
 \endproclaim\startproof
 We treat only the case where $L/K$ is a finite extension.  (The
 general case can surely be reduced to this one.)
 Let $X$ be the abelian variety~$X=\Res_{L/K} A$.
 As in~\S6, we find a decomposition
 $$ \Q\otimes\End_K(X) = \prod_{\sigma\in\Gal(L/K)}
 \Q\otimes\Hom_L(\ssigma A,A).$$
 The homomorphism $\mu_\sigma$ in
 the ``$\sigma^{\roman{th}}$ factor" corresponds to an
 element $[\sigma]$ of~$\Q\otimes\End_K(X)$.
 This element operates on $\prod \ssigma A$ as
 a matrix, sending ${}^{\tau g}\mkern-\thinmuskip A$ to
 $\ttau A$ by $\ttau\mu_g$ for each $\tau\in\Gal(L/K)$.

 Because
 the analogue of~$c(\sigma,\tau)$ is~1 in this context, we have
 simply $[\sigma][\tau]=[\sigma\tau]$ for $\sigma,\tau\in\Gal(L/K)$:
 the identity~(8.1) shows that the product of the
 matrices representing $[\sigma]$ and~$[\tau]$ is the matrix
 representing $[\sigma\tau]$.
 Hence $\Cal R\seteq\Q[\Gal(L/K)]$ operates on~$X$.

 An analogue of Proposition~6.5 shows that we have
 $$ X_L \approx \Cal R\otimes_\Q A $$
 in the category of abelian varieties over~$L$ up to isogeny.
 To see this explicitly, we let $A_\sigma$ be a copy of~$A$
 indexed by $\sigma$ and consider the isomorphism of abelian
 varieties over~$L$ up to isogeny
 $$\iota\: \prod_\sigma A_\sigma \,\buildrel\sim\over\to\,
       \prod_\sigma \ssigma A = X_L$$
 which takes $A_\sigma$ to ${}^{\sigma^{-1}}\mkern-\thinmuskip A$
 via the map ${}^{\sigma^{-1}}\mkern-\thinmuskip \mu_\sigma$.
 Via this isomorphism, the automorphism $[\sigma]$ of~$X$ is
 transported to the permutation which takes each factor $A_g$ of
 $\prod_g A_g$ to the factor $A_{\sigma g}$, via the identity map
 $A\to A$.

 Let $B$ be the image of $\eta\seteq \sum_\sigma [\sigma]$, so that
 $B$ is an abelian subvariety of~$X$ which is defined
 over~$K$.  (The sum $\eta$ need not
 be a literal endomorphism of~$X$, so that, strictly speaking, one
 should consider the image of a suitable multiple of~$\eta$.)
 The isomorphism $\iota$ sends $B_L$ to the diagonal image of~$A$
 in $\prod A_\sigma = A\times\cdots\times A$.  Hence $B_L$ is isogenous
 to~$A$.
 \endproof

 \enddocument